\documentclass[11pt]{article}
\usepackage{amssymb}

\begin{document}

\newcommand{\Rm}{{\mathbb R}}
\newcommand{\disp}{\displaystyle}
\newcommand{\eps}{\varepsilon}

\title{Global Well-posedness and Scattering of Defocusing Energy subcritical Nonlinear Wave Equation in dimension 3 with radial data}
\author{Ruipeng Shen,\\
        Department of Mathematics,\\
        University of Chicago}
\date{November 2011}
\maketitle
\section{Introduction}
In this paper we will consider the defocusing case of the energy subcritical non-linear wave equation in $\Rm^3$ with radial initial data.
\begin{equation}
\left\{\begin{array}{l} \partial_t^2 u - \Delta u + |u|^p u = 0, \,\,\,\, (x,t)\in \Rm^3 \times \Rm,\\
u(0) = u_0 \in \dot{H}^s (\Rm^3),\\
\partial_t u(0) = u_1 \in \dot{H}^{s-1}(\Rm^3).\end{array}\right. 
\label{eqn}
\end{equation}
Here 
\[
 p = \frac{2}{3/2 - s}.
\]
The case $s=1$ is the energy-critical case. The following quantity is called the energy of the solution. 
The energy is constant for all time, as long as the solution still exists. 
\begin{equation}
 E(t) = \frac{1}{2} \int_{\Rm^3} \left(|\partial_t u (x,t)|^2 + |\nabla u (x,t)|^2\right) dx + 
\frac{1}{6} \int_{\Rm^3} |u(x,t)|^6 dx. \label{defenergy}  
\end{equation}
Thanks to the existence of the energy, one can show the universal boundedness of the following norms for all time with ease.
\[
 \|\partial_t u\|_{L^2}; \|\nabla u\|_{L^2}; \|u\|_{L^6}.
\]
The global well-posedness and scattering of the solutions in the energy-critical case is known. Please see \cite{mg1,mg2}.\\
In this paper we will consider the case when $s$ is slightly smaller than $1$. In this case
the energy does not exist. Thus we can not obtain the similar boundedness as the energy critical case. Instead we have 
to make the following assumption 
\begin{equation}
 \sup_{t \in I} \|(u, \partial_t u)\|_{{\dot{H}^s}\times{\dot{H}^{s-1}}} < \infty. \label{add1}
\end{equation}
where $I$ is the maximal interval of existence of the solution.
\paragraph{Remark} Please note that although this problem is energy subcritical, it is actually 
$\dot{H}^s \times \dot{H}^{s-1}$ critical by the choice of $p$, because if $u(x,t)$ is a solution
of (\ref{eqn}) with initial data $(u_0,u_1)$, then for any $\lambda>0$,
\[
 \frac{1}{\lambda^{3/2 - s}} u (\frac{x}{\lambda}, \frac{t}{\lambda})
\]
is another solution of the equation (\ref{eqn}) with the initial data
\[
 \left(\frac{1}{\lambda^{3/2 - s}} u_0 (\frac{x}{\lambda}), 
\frac{1}{\lambda^{5/2 - s}} u_1 (\frac{x}{\lambda})\right).
\]
These two pairs of initial data share the same $\dot{H}^s \times \dot{H}^{s-1}$ norm.
These scalings play an important role in our discussion of this problem.
\paragraph{Main theorem} Assume $s > 15/16$. Let $u(t)$ be a solution of (\ref{eqn}) with a maximal life span $I$ and radial
initial data $(u_0,u_1) \in \dot{H}^s \times \dot{H}^{s-1}$. In addition, we assume $u(t)$ satisfies the uniform boundedness condition
(\ref{add1}). Then $u(t)$ is a global solution (i.e $I = \Rm$) and scatters. 
\[
 \|u\|_{L^{2p} L^{2p} (\Rm \times \Rm^3)} < \infty.
\] 
This is actually equivalent to saying that there exist two pairs $(u_0^+,u_1^+)$ and $(u_0^-,u_1^-)$ 
in $\dot{H}^s \times \dot{H}^{s-1}$ such that
\[
 \lim_{t \rightarrow \pm \infty} \|u(t) - S(t)(u_0^\pm, u_1^\pm)\|_{\dot{H}^s \times \dot{H}^{s-1}} = 0.
\]
Here $S(t)(u_0,u_1)$ is the solution of the Linear Wave Equation with the initial data $(u_0, u_1)$.

\paragraph{The structure of the paper} We will introduce some a priori estimates and then introduce the local theory at the very 
beginning. The main idea to prove the main scattering result is to show:\\
(I) If the theorem failed, it would break down for a special solution with a critical norm.\\
(II) The solution in (I) does not exist.\\
The first step in this process is somewhat standard to deal with similar problems of dispersive equations. 
Thus we will only give important statements instead of showing all the details in this step. One could 
read the references if one is interested in how to establish these results.\\
The second step, however, depends on the specific problems. Thus the majority of this paper consists of 
concrete discussions of this step.

\section{A Priori Results}
In this section, we will review the theory for the Cauchy problem of nonlinear wave equation (\ref{eqn}).\\
Let $I$ be an interval of time. We define the following norms with $\frac{1}{2}\leq s \leq 1$  
\[
 \|v\|_{S(I)} = \|v\|_{L^{2p} L^{2p} (I \times \Rm^3)};
\]
\[
 \|v\|_{W(I)} = \|v\|_{{L^4} L^{4} (I \times \Rm^3)}.
\] 
The space-time norm is defined by
\[
 \|v(x,t)\|_{L^{q} L^{r} (I \times \Rm^3)} = \left(\int_I \left(\int_{\Rm^3} |v(x,t)|^r dx\right)^{q/r} dt\right)^{1/q}.
\]
We say $u(t) (t\in I)$ is a solution of (\ref{eqn}), if $(u,\partial_t u) \in C(I;{\dot{H}^s}\times{\dot{H}^{s-1}} )$, with finite norms $\|u\|_{S(J)}$ and $\|D_x^{s - {1/2}} u\|_{W(J)}$ for any bounded closed interval $J \subseteq I$ so that the integral equation 
\[
 u(t) = S(t)(u_0,u_1) + \int_0^t \frac{\sin ((t-s)\sqrt{-\Delta})}{\sqrt{-\Delta}} F(u(s)) ds
\]
holds for all time $t \in I$.  
\[
 F(u) = -|u|^{p}u.
\]
\paragraph{Generalized Strichartz Inequalities}. (Please see proposition 3.1 of \cite{strichartz}, here we use the Sobolev version in $\Rm^3$)
Let $2 \leq q_1,q_2 \leq \infty$, $2 \leq r_1, r_2 < \infty$ and $\rho_1, \rho_2, s \in \Rm$ with
\[
 1/{q_i} + 1/{r_i} \leq 1/2; \,\,\,\, i=1,2.
\]
\[
 1/{q_1} + 3/{r_1} = 3/2 - s + \rho_1.
\]
\[
 1/{q_2} + 3/{r_2} = 1/2 + s + \rho_2. 
\]
In particular, if $(q_1, r_1, s, \rho_1) = (q, r, m, 0)$ satisfies the conditions above, we say
$(q,r)$ is an $m$-admissible pair.\\
Let $u$ be the solution of the following linear wave equation 
\begin{equation}
\left\{\begin{array}{l} \partial_t^2 u - \Delta u = F (x,t), \,\,\,\,\, (x,t)\in \Rm^3 \times \Rm;\\
u(0) = u_0 \in \dot{H}^s (\Rm^3);\\
\partial_t u(0) = u_1 \in \dot{H}^{s-1}(\Rm^3).\end{array}\right. 
\end{equation}
Then we have 
\begin{eqnarray*}
 && \|(u(T), \partial_t u(T))\|_{\displaystyle \dot{H}^s \times \dot{H}^{s-1}} + \|D_x^{\rho_1} u\|_{\displaystyle L^{q_1} L^{r_1} ([0,T]\times \Rm^3)}\\
 &\leq& 
 C \left( \|(u_0,u_1)\|_{\displaystyle \dot{H}^s \times \dot{H}^{s-1}} + \|D_x^{-\rho_2} F(x,s)\|_{\displaystyle L^{\bar{q}_2} L^{\bar{r}_2}([0,T]\times \Rm^3)} \right). 
\end{eqnarray*}
The constant $C$ does not depend on $T$.\\
Using the Strichartz estimate and a fixed-point argument, we have the following theorems. 
(Please see \cite{ls} for more details) 
\paragraph{Theorem 1(Local solution)} For any initial data $(u_0,u_1) \in \dot{H}^s \times \dot{H}^{s-1}$, there is a maximal interval $(-T_{-}(u_0,u_1), T_{+}(u_0,u_1))$ in which the equation has a solution. 
\paragraph{Theorem 2(Scattering with small data)} There exists $\delta > 0$ such that if $\|(u_0,u_1)\|_{\dot{H}^s \times \dot{H}^{s-1}} < \delta$, then the Cauchy problem (\ref{eqn}) has a global-in-time solution $u$ with $\|u\|_{S(-\infty,+\infty)} < \infty$. 
\paragraph{Lemma(Standard finite blow-up criterion)} If $T_{+} < \infty$ under the uniform boundedness condition (\ref{add1}), then 
\[
 \|u\|_{S([0,T_{+}))} = \infty.
\]
\paragraph{Theorem 3(Long time perturbation theory)}(See \cite{ctao, kenig, kenig1, km}) Let $M, A, A'$ be positive constants. There exists $\eps_0 = \eps_0(M,A,A')>0$ and $\beta>0$ such that if $\eps < \eps_0$, for any approximation solution $\tilde{u}$ defined on $\Rm^3 \times I$
and any initial data $(u_0,u_1) \in \dot{H}^s \times \dot{H}^{s-1}$ satisfying
\[
 ({\partial_t}^2 - \Delta) (\tilde{u}) - F(\tilde{u}) = e, \,\,\,\,\, (x,t) \in \Rm^3 \times I; 
\]
\begin{equation}
\left\{ \begin{array}{l}
\sup_{t \in I} \|(\tilde{u}(t), \partial_t \tilde{u}(t))\|_{\dot{H}^s \times \dot{H}^{s-1}} \leq A,\\
 \|\tilde{u}\|_{S(I)} \leq M,\\
\|\tilde{u}\|_{W(J)} < \infty \hbox{ for each } J\subset\subset I; \end{array} \right. 
\end{equation}
\[
 \|(u_0-\tilde{u}(0), u_1-\partial_t \tilde{u}(0))\|_{\dot{H}^s \times \dot{H}^{s-1}} \leq A';
\]
\[
 \|D_x^{s - \frac{1}{2}} e\|_{L_I^{4/3} L_x^{4/3}} + \|S(t)(u_0-\tilde{u}(0),u_1 - \partial_t \tilde{u}(0))\|_{S(I)} \leq \eps.
\]
Then there exists a solution of (\ref{eqn}) defined in the interval $I$ with the initial date $(u_0,u_1)$ and satisfying
\[
  \|u\|_{S(I)} \leq C(M,A,A');
\]
\[
 \sup_{t\in I} \|((u(t), \partial_t u(t)) - ((\tilde{u}(t), \partial_t \tilde{u}(t))\|_{\dot{H}^s \times \dot{H}^{s-1}} \leq C(M,A,A')(A' +\eps +\eps^\beta).
\]
\paragraph{Remark} If $K$ is a compact set in the space $\dot{H}^s \times \dot{H}^{s-1}$, then there exists $T = T(K) > 0$ such that for any $(u_0,u_1) \in K$, $T_{+}(u_0,u_1) > T(K)$. This
is a direct result from the perturbation theory.

\paragraph{A Global Integral Estimate} At the end of this section we have a global integral estimate for the solution $u$. Unlike the local theory, this estimate could only be applied to
a solution in the energy space. 
\paragraph{Lemma} (Please see \cite{benoit}) Let $u$ be a solution of (\ref{eqn}) defined in a time interval $[0,T]$ with $(u,\partial_t u) \in \dot{H}^1 \times L^2$ and a finite energy
\[
 E = \int_{\Rm^3} (|\nabla_x u|^2 + |\partial_t u|^2 + \frac{2}{p+2}|u(x)|^{p+2} ) dx.
\]
For any $R>0$, we have 
\begin{eqnarray*}
 \frac{1}{2R} \int_0^T \int_{|x|< R} (|\nabla u |^2 + |\partial_t u|^2) dx dt + \frac{1}{2R^2} \int_0^T \int_{|x|=R} |u|^2 d\sigma_R dt && \\
+ \frac{1}{2R} \frac{2p - 2}{p+2} \int_0^T \int_{|x|<R} |u|^{p+2} dx dt
+ \frac{p}{p+2} \int_0^T \int_{|x|>R} \frac{|u|^{p+2}}{|x|} dx dt &&\\ + \frac{2}{R^2} \int_{|x|<R} |u(T)|^2
 &\leq& E.
\end{eqnarray*}
Observing that each term on the left hand is nonnegative, we can obtain a uniform upper bound for the last term in the second line above
\[
  \int_0^T \int_{|x|>R} \frac{|u|^{p+2}}{|x|} dx dt \leq \frac{p+2}{p} E.
\]
Let $R$ approach zero and $T$ approach $T_{+}$, we have 
\begin{equation}
\int_0^{T_{+}} \int_{\Rm^3} \frac{|u|^{p+2}}{|x|} dx dt \leq \frac{p+2}{p} E.
\label{upperbound}
\end{equation}
\section{Compactness Process}
As we stated in the first section, the standard technique here is to show if the main theorem failed, there would be a special minimal blow-up solution. 
In addition, this solution is almost periodic modulo symmetries. Namely the set 
\[
 \left\{\left(\frac{1}{\lambda(t)^{3/2-s}} u(\frac{x}{\lambda(t)}, t), \frac{1}{\lambda(t)^{5/2-s}} \partial_t u(\frac{x}{\lambda(t)}, t)\right): t \in I\right\}
\]
is precompact in $\dot{H}^s \times \dot{H}^{s-1}$. The function $\lambda(t)$ is called the frequency scale function,
because the solution $u(t)$ at time $t$ concentrates around the frequency $\lambda(t)$.\\ 
Please note that here we use 
the radial condition, thus the only available symmetries are scalings. If we did not assume the radial condition, similar
results would still hold but the symmetries would include translations besides scalings.\\
  
The following is the first compactness result. 
\paragraph{Minimal blow-up solution} Assume that the main theorem failed. Then there would exist
a solution $u: I \times \Rm^3 \rightarrow \Rm$ such that
\[
\sup_{t\in I} \|(u, \partial_t u)\|_{\dot{H}^{s} \times \dot{H}^{s-1}} < \infty,
\]
$u$ blows up in the positive direction at time $T_+ \leq + \infty$ with
\[
 \|u\|_{S(0, T_+)} = \infty.
\]
In addition, $u$ is almost periodic modulo scaling with a frequency scale function $\lambda(t)$.
It is minimal in the following sense, if 
\[
 \sup_{t \in J} \|(v, \partial_t v)\|_{\dot{H}^{s} \times \dot{H}^{s-1}} < \sup_{t\in I} \|(u, \partial_t u)\|_{\dot{H}^{s} \times \dot{H}^{s-1}}
\]
for another solution $v$ with a maximal lifespan $J$, then $v$ is a global solution in time  and scatters. \\
The main tool to obtain this result is the profile decomposition. One could follow the argument in \cite{kenig2} in order to find a proof. 
In that paper C.E.Kenig and F.Merle deal with the cubic defocusing NLS under similar assumptions.\\
The second compactness result is that we can always assume the frequency scale function $\lambda(t)$ is uniformly bounded for all $t \geq 0$.
If this was not true for our minimal blow-up solution $u$ mentioned above, one could always take a sequence
\[
 0 < t_1 < t_2 < t_3 < \cdots < t_n \cdots <  T_+,
\]
such that $\lambda(t) \leq \lambda(t_i)$ for all $0 \leq t \leq t_i$ and 
\[
 \lim_{i \rightarrow \infty} \lambda(t_i) = \infty.
\]
Using the compactness results, we know
\[
 \left\{\left(\frac{1}{\lambda(t_i)^{3/2-s}} u(\frac{x}{\lambda(t_i)}, t_i), -\frac{1}{\lambda(t_i)^{5/2-s}} \partial_t u(\frac{x}{\lambda(t_i)}, t_i)\right): i \in {\mathbb Z}^+\right\}
\]
is precompact in $\dot{H}^{s} \times \dot{H}^{s-1}$. Thus one can extract a subsequence so that
this subsequence converges to a pair $(v_0,v_1)$ in $\dot{H}^{s} \times \dot{H}^{s-1}$.\\
Consider the solution $v$ of the original equation with the initial data $(v_0,v_1)$, we can show $v$ is 
still a minimal blow-up solution and its frequency scale function $\lambda(t) \leq 1$ for all $t\geq 0$. (Please
see \cite{kenig,kenig1} for more details)\\
Using the remark following the perturbation theory, one can obtain $T_+ = \infty$ for our 
minimal solution $v$ immediately from the fact that $\lambda(t) \leq 1$. \\
In summary, if the main theorem failed, we would find a minimal blow-up solution $u$, so that 
it blows up at $T_+ = \infty$ and its frequency scale function $\lambda(t) \leq 1$ for all $t \geq 0$.


 
\paragraph{Local Compactness} 
Let $s \geq 3/4$. Fix a cutoff function $\varphi(x) \in C^{\infty}$ with the following properties.\\
\[
 \varphi(x) \left\{\begin{array}{l} =0, \,\,\,\, |x| \leq 1/2;\\
 \in [0,1],\,\,\,\, 1/2 \leq |x| \leq 1;\\
 = 1,\,\,\,\, |x| \geq 1.\end{array}\right. 
\]
For a minimal blow-up solution mentioned above and its frequency scale function $\lambda(t)$, we have the following 
propositions by a compactness argument. 
\paragraph{Proposition} If $u\neq 0$, there exist $d, C', R_1 > 0$ and $C_1 > 1$ independent of $t$ such that\\
(i) The interval $[t-d \lambda^{-1}(t), t+ d \lambda^{-1}(t)] \subseteq I$ for all $t \in I$. In addition, for all $s \in [t-d \lambda^{-1}(t), t+ d \lambda^{-1}(t)]$,
\begin{equation}
 \frac{1}{C_1} \lambda(t) \leq \lambda(s) \leq C_1 \lambda(t).
\label{choiced}
\end{equation}
(ii) We have the following estimate for an $s$-admissible pair $(q,r)$.
\[
 \|u\|_{\displaystyle L^q  L^{r} ([t-d \lambda^{-1}(t), t+ d \lambda^{-1}(t)] \times \Rm^3 )} \leq C'.
\]
(iii) 
\[
 \left\|\left(\varphi(\frac{x}{R_1 \lambda^{-1}(t)})u,\varphi(\frac{x}{R_1 \lambda^{-1}(t)})\partial_t u\right)\right\|_{\dot{H}^1 \times \dot{H}^{s-1}} \leq \delta.  
\]
Here $\delta$ is the small constant we need to apply the global solution theory for small data. \\
(iv) We can also get a lower bound. By a compactness argument we obtain that there exist $R_0, \eta_0>0$, so that for all $t$,
\[
 \int_0^d \int_{|x| < R_0} \frac{(\frac{1}{\lambda(t)^{2/p}} |u(\lambda^{-1}(t)x, \lambda^{-1}(t)s + t)|)^{p+2}}{|x|}dxds \geq \eta_0. 
\]
This implies
\[
 \int_0^d \int_{|x| < R_0} \frac{ |u(\lambda^{-1}(t)x, \lambda^{-1}(t)s + t)|^{p+2}}{\lambda^{-1}(t)|x|} \frac{dxds}{\lambda(t)^{\frac{2}{p}(p+2)+1}} \geq \eta_0. 
\]
\[
 \frac{1}{\lambda(t)^{4/p -1}} \int_0^d \int_{|x| < R_0} \frac{ |u(\lambda^{-1}(t)x, \lambda^{-1}(t)s + t)|^{p+2}}{\lambda^{-1}(t)|x|} \frac{dxds}{\lambda(t)^4} \geq \eta_0. 
\]
\begin{equation}
\begin{array}{r}
  \displaystyle \int_{t}^{t + d \lambda^{-1}(t)} \int_{|x| < R_0 \lambda^{-1}(t)} \frac{|u(x,s)|^{p+2}}{|x|} dx ds \geq  \lambda(t)^{4/p -1} \eta_0\\
   =  \lambda(t)^{2 -2s} \eta_0.
\end{array} \label{lowbound1}
\end{equation}
This lower bound is essential to give a contradiction in the last part of this paper. Please see (\cite{kenig2}) for more details
of this kind of argument.

\section{Regularity of Solutions} In this section, we will show the solution $u$ we obtained in the previous section has additional regularity by the following Duhamel formula.
The additional regularity will enable us to use the methods and estimates only available in the energy space.
\begin{equation}
 u(t) = \int_t^{+\infty} \frac{\sin((s-t)\sqrt{-\Delta})}{\sqrt{-\Delta}} F(u(s)) ds. 
\label{duhamel}
\end{equation}
\[
 \partial_t u(t) = -\int_t^{+\infty} \cos((s-t)\sqrt{-\Delta}) F(u(s)) ds.
\]
These identities hold in the sense of weak limits in $\dot{H}^{s} \times \dot{H}^{s-1}$. Please note that the minimal solution we obtained in the previous section exists forever in the positive direction. 
This fact makes it possible for us to take the integral from $t$ to $\infty$. In addition,
given fixed closed interval $J$ compactly 
supported in the maximal lifespan $I$, the first identity also holds in the sense of strong limit in the space $L^q L^r (J \times \Rm^3)$ for an admissible pair $(q,r)$ with $q < \infty$.\\
\subsection{Local Contribution} 
\paragraph{New Norm} Let us define the following $X(J)$ norm for an interval $J$ contained in the maximal lifespan.
\[
 \|u\|_{X(J)} = \|u\|_{\displaystyle L^{p+1} L^{\frac{6(p+1)}{5-2s}} (J \times \Rm^3)}. 
\]
If $s=1$, this is the classic $L^5 L^{10}$ norm. This pair is admissible as long as  
\[
 s > \frac{11 - \sqrt{73}}{4}. 
\]
At this time, let us choose $s > 3/4$, so the $X(J)$ norm can be estimated by the Strichartz Inequality.
\paragraph{Definition} Let us define
\[
 M(A) = \sup_{t \in I} \|(u_{> \lambda(t)A}, P_{> \lambda(t)A} \partial_t u)\|_{\displaystyle \dot{H}^s \times \dot{H}^{s-1}}.
\]
\[
 S(A) = \sup_{t \in I} \|u_{> \lambda(t)A}\|_{\displaystyle X([t, t+ d \lambda^{-1}(t)])}.
\]
\[
 N(A) = \sup_{t \in I} \|P_{> \lambda(t)A}(F(u))\|_{\displaystyle L^1 L^{\frac{6}{5-2s}}([t, t+ d\lambda^{-1}(t)]\times \Rm^3)}.
\]
The Operator $P$ is the smooth frequency cutoff operator. While the subscript of $u$ has the same meaning.
\[
 u_{> \lambda(t)A} = P_{> \lambda(t)A} u.
\] 
By a compactness argument, these $N(A), M(A), S(A)$ are bounded by a universal constant for all $A>0$. 
They tend to $0$ as $A$ goes to $\infty$. Our goal is to gain decay of $S(A)$ and $N(A)$.
This decay will give us some additional regularity of the solution.\\ 
First we need to prove two technical lemmas used in the argument. 
\paragraph{Lemma} Let $u$ be a function defined in $\Rm^3 \times J$ and $s > 3/4$. Suppose the support of $\hat{u}$ is contained in the ball $B(0,r)$
for each $t \in J$. Then 
\[
 \|P_{>R} F(u)\|_{L^1 L^{\frac{6}{5-2s}}(J \times \Rm^3) } \lesssim (\frac{r}{R})^{2-s} \|u\|_{L^p L^{3p}(J \times \Rm^3)}^p \sup_{t \in J} \|u\|_{\dot{H}^s}.
\]
\paragraph{Proof}: Using the Sobolev embedding and the frequency cutoff we can estimate the left hand by
\begin{eqnarray*}
 \lefteqn{\| P_{>R} (|u|^p u)\|_{L^1 L^{\frac{6}{5-2s}}}
 \lesssim \|D^s_x P_{>R}(|u|^p u)\|_{L^1 L^{6/5}}}\\
&\lesssim& \frac{1}{R^{2 -s}} \|\Delta_x P_{>R}(|u|^p u)\|_{L^1 L^{6/5}}\\
&\lesssim& \frac{1}{R^{2 -s}} \|\Delta_x (|u|^p u)\|_{L^1 L^{6/5}}\\
&\lesssim& \frac{1}{R^{2 -s}} \|(|u|^{p-2}u)|\nabla u|^2\|_{L^1 L^{6/5}}
  + \frac{1}{R^{2 -s}} \||u|^p (\Delta_x u)\|_{L^1 L^{6/5}}.
\label{ls1}
\end{eqnarray*}
By the frequency cutoff and the Bernstein Inequality, we have 
\[
 \|\nabla u\|_{L^{\infty} L^2} \lesssim r^{1-s} \|D_x^s u\|_{L^{\infty} L^2} \leq r^{1-s} \sup_{t} \|u\|_{\dot{H}^s}.
\]
\[
 \|\nabla u\|_{L^p L^{3p}} \lesssim r \|u\|_{L^p L^{3p}}.
\]
Thus the first term can be estimated by 
\begin{eqnarray*}
 \|(|u|^{p-2}u)|\nabla u|^2\|_{L^1 L^{6/5}} &\lesssim& \|u\|_{L^p L^{3p}}^{p-1} \|\nabla u\|_{L^p L^{3p}} \|\nabla u\|_{L^\infty L^2}\\
&\lesssim&  \|u\|_{L^p L^{3p}}^{p-1} ( r \|u\|_{L^p L^{3p}}) ( r^{1-s} \sup_{t} \|u\|_{\dot{H}^s})\\
&\lesssim& r^{2-s} \|u\|_{L^p L^{3p}}^{p} \sup_{t} \|u\|_{\dot{H}^s}.
\end{eqnarray*}
We also have 
\[
 \|\Delta_x u\|_{L^\infty L^2} \lesssim r^{2-s} \|D_x^s u\|_{L^\infty L^2} = r^{2-s} \sup_{t} \|u\|_{\dot{H}^s}.
\]
This gives the estimate of the second term 
\begin{eqnarray*}
 \||u|^p (\Delta_x u)\|_{L^1 L^{6/5}} &\lesssim& \|u\|_{L^p L^{3p}}^p \|\Delta_x u\|_{L^\infty L^2}\\
&\lesssim& r^{2-s} \|u\|_{L^p L^{3p}}^p \sup_{t} \|u\|_{\dot{H}^s}.
\end{eqnarray*}
Combining these estimates, we have the inequality in the lemma. 
\paragraph{Bilinear Estimate} Suppose $u_i$ satisfy the following linear wave equation on the time interval $I = [0,T]$
\[
 \partial_{tt} u_i - \Delta u_i = F_i (x,t),
\]
with the initial data $(u_i (0), \partial_t u_i (0)) = (u_{0,i}, u_{1,i})$. Then 
\begin{eqnarray*}
 S &=& \|(P_{>R} u_1) (P_{<r}u_2)\|_{L^{\frac{p+1}{2}} L^{\frac{3(p+1)}{5 -2s}} (I \times \Rm^3)} \\
 & \lesssim & (\frac{r}{R})^\sigma \left(\|(u_{0,1}, u_{1,1})\|_{\dot{H}^s \times \dot{H}^{s-1}} + \|F_1\|_{L^1 L^{\frac{6}{5-2s}}(I \times \Rm^3)}\right)\\ 
&& \times \left(\|(u_{0,2}, u_{1,2})\|_{\dot{H}^s \times \dot{H}^{s-1}} + \|F_2\|_{L^1 L^{\frac{6}{5-2s}}(I \times \Rm^3)}\right).
\end{eqnarray*}
This estimate is meaningful only if $r\ll R$. Otherwise it is just the Strichartz estimate. 
Here the number $\sigma$ is given by the following
\begin{equation}
 0 < \sigma < 3 \min \left\{\frac{1}{2} - \frac{1}{p+1} - \frac{5-2s}{6(p+1)}, \frac{5 -2s}{6(p+1)}\right\}.
\label{sigma}
\end{equation}
\paragraph{Proof} By the Strichartz estimate 
\begin{eqnarray*}
 &&\|(P_{>R}) u_1\|_{\displaystyle L^{p+1} L^{1/{(\frac{5-2s}{6(p+1)} + \frac{\sigma}{3}})}}\\
 &\lesssim& \|(D_x^{-\sigma} P_{>R} u_{0,1}, D_x^{-\sigma} P_{>R} u_{1,1})\|_{\dot{H}^s \times \dot{H}^{s-1}} + \|D_x^{-\sigma} P_{>R} F_1\|_{L^1 L^{\frac{6}{5-2s}}}.
\end{eqnarray*}
\begin{eqnarray*}
 &&\|(P_{<r}) u_2\|_{\displaystyle L^{p+1} L^{1/{(\frac{5-2s}{6(p+1)} - \frac{\sigma}{3}})}}\\
 &\lesssim& \|(D_x^{\sigma} P_{<r} u_{0,2}, D_x^{\sigma} P_{<r} u_{1,2})\|_{\dot{H}^s \times \dot{H}^{s-1}} + \|D_x^{\sigma} P_{<r} F_2\|_{L^1 L^{\frac{6}{5-2s}}}.
\end{eqnarray*}
Our choice of $\sigma$ makes sure that the pairs above are admissible. Thus we have  
\begin{eqnarray*}
 && \|(P_{>R} u_1) (P_{<r}u_2)\|_{\displaystyle L^{\frac{p+1}{2}} L^{\frac{3(p+1)}{5 -2s}}} \\
&\lesssim& \|(P_{>R}) u_1\|_{\displaystyle L^{p+1} L^{1/{(\frac{5-2s}{6(p+1)} + \frac{\sigma}{3})}}} \|(P_{<r}) u_2\|_{\displaystyle L^{p+1} L^{1/{(\frac{5-2s}{6(p+1)} - \frac{\sigma}{3})}}}\\
&\lesssim& \left(\|(D_x^{-\sigma} P_{>R} (u_{0,1}), D_x^{-\sigma} P_{>R} (u_{1,1}))\|_{\dot{H}^s \times \dot{H}^{s-1}} + \|D_x^{-\sigma} P_{>R} F_1\|_{L^1 L^{\frac{6}{5-2s}}}\right)\\
&& \times \left(\|(D_x^{\sigma} P_{<r} (u_{0,2}), D_x^{\sigma} P_{<r} (u_{1,2}))\|_{\dot{H}^s \times \dot{H}^{s-1}} + \|D_x^{\sigma} P_{<r} F_2\|_{L^1 L^{\frac{6}{5-2s}}}\right)\\
&\lesssim& (\frac{1}{R})^{\sigma} \left(\|( P_{>R} (u_{0,1}), P_{>R} (u_{1,1}))\|_{\dot{H}^s \times \dot{H}^{s-1}} + \| P_{>R} F_1\|_{L^1 L^{\frac{6}{5-2s}}}\right)\\
&& \times r^\sigma \left(\|( P_{<r} (u_{0,2}), P_{<r} (u_{1,2}))\|_{\dot{H}^s \times \dot{H}^{s-1}} + \|P_{<r} F_2\|_{L^1 L^{\frac{6}{5-2s}}}\right)\\
& \lesssim & \hbox{the right hand}.
\end{eqnarray*}
\paragraph{Recurrence Formulas} Let $s >3/4$. We have the following formulas for any $0 < \alpha < \beta <1$, positive constant $\eps_1$ and sufficiently large $A$.
These formulas are essential in our proof of the decay.
\begin{equation}
 N(A) \lesssim S(A^\beta) S^{p}(A^\alpha) + A^{-(\beta - \alpha)\sigma} + A^{-(2-s)(1-\beta)}.
 \label{re1}
\end{equation}
\begin{equation}
 S(A) \lesssim N(A^{1-\eps_1}) + A^{-\sigma_1}.
\label{re2}
\end{equation}
The constant $\sigma$ is the one in our bilinear estimate. While the constant $\sigma_1$ in the second inequality is given by
\[
 \sigma_1 = \frac{4(sp +s-2)}{3p(p+1)} > 0. 
\]
\paragraph{Proof of the First Inequality} To prove formula (\ref{re1}), we chop the solution in frequency. All the norms below are taken in the time-space $[t,t+ d\lambda^{-1}(t)]\times \Rm^3$.
\begin{eqnarray*}
 \|P_{> \lambda(t)A}(F(u))\|_{L^1 L^{\frac{6}{5-2s}}} &\leq& \|P_{> \lambda(t)A}(F(u) - F(u_{\leq A^\beta \lambda(t)}))\|_{L^1 L^{\frac{6}{5-2s}}}\\
 && +  \|P_{> \lambda(t)A}(F(u_{\leq A^\beta \lambda(t)}))\|_{L^1 L^{\frac{6}{5-2s}}}.
\end{eqnarray*}
The first term equals
\begin{eqnarray*}
 &&\left\|P_{>\lambda(t)A}[u_{>A^\beta \lambda(t)} \int_0^1 F'(u_{\leq A^\beta \lambda(t)} + s u_{> A^\beta \lambda(t)})ds]\right\|_{L^1 L^{\frac{6}{5-2s}}}\\
&\lesssim& \left\|u_{>A^\beta \lambda(t)} \int_0^1 F'(u_{\leq A^\beta \lambda(t)} + s u_{> A^\beta \lambda(t)})ds\right\|_{L^1 L^{\frac{6}{5-2s}}}\\
&\lesssim& \left\| \begin{array}{l} u_{>A^\beta \lambda(t)} \displaystyle \int_0^1 F'(u_{\leq A^\beta \lambda(t)} + s u_{> A^\beta \lambda(t)}) ds\\
 -u_{>A^\beta \lambda(t)} \displaystyle \int_0^1 F'(u_{ A^\alpha \lambda(t)< \cdotp\leq A^\beta \lambda(t)} + s u_{> A^\beta \lambda(t)}) ds \end{array} \right\|_{L^1 L^{\frac{6}{5-2s}}}\\
&& + \left\|u_{>A^\beta \lambda(t)} \int_0^1 F'(u_{ A^\alpha \lambda(t)<\cdotp\leq A^\beta \lambda(t)} + s u_{> A^\beta \lambda(t)})ds\right\|_{L^1 L^{\frac{6}{5-2s}}}\\
&\lesssim& \left\| \begin{array}{l} u_{>A^\beta \lambda(t)} u_{\leq A^\alpha \lambda(t)} \\
 \times \displaystyle \int_0^1 \int_0^1  F''( \tilde{s} u_{\leq A^\alpha \lambda(t)} + u_{ A^\alpha \lambda(t)< \cdotp\leq A^\beta \lambda(t)} + s u_{> A^\beta \lambda(t)} ) ds d\tilde{s} \end{array} \right\|\\
&& + \left\|u_{>A^\beta \lambda(t)} \int_0^1 F'(u_{ A^\alpha \lambda(t)<\cdotp\leq A^\beta \lambda(t)} + s u_{> A^\beta \lambda(t)})ds\right\|_{L^1 L^{\frac{6}{5-2s}}}\\
&\lesssim& \left\|u_{>A^\beta \lambda(t)} u_{\leq A^\alpha \lambda(t)}\right\|_{\displaystyle L^{\frac{p+1}{2}} L^{\frac{3(p+1)}{5-2s}}}\\
&& \times \left\|\begin{array}{r} \displaystyle \int_0^1 \int_0^1 F''(\tilde{s} u_{\leq A^\alpha \lambda(t)} + u_{ A^\alpha \lambda(t)< \cdotp\leq A^\beta \lambda(t)}\\
+ s u_{> A^\beta \lambda(t)}) ds d\tilde{s}\end{array} \right\|_{\displaystyle L^{\frac{p+1}{p-1}} L^{ \frac{6(p+1)}{(p-1)(5-2s)}}}\\
&+& \left\| u_{>A^\beta \lambda(t)} \right\|_{L^{p+1} L^{\frac{6(p+1)}{5-2s}}}\\
&& \times \left\|\int_0^1 F'(u_{ A^\alpha \lambda(t)<\cdotp\leq A^\beta \lambda(t)} + s u_{> A^\beta \lambda(t)})ds\right \|_{L^{\frac{p+1}{p}} L^{\frac{6(p+1)}{p(5-2s)}}}\\
&\lesssim& (\frac{A^\alpha \lambda(t)}{A^\beta \lambda(t)})^\sigma + S(A^\beta) S^p (A^\alpha)\\
&\lesssim& A^{-(\beta-\alpha)\sigma} + S(A^\beta)S^p(A^\alpha).
\end{eqnarray*}
The bilinear estimate is used here in order to estimate the term $u_{>A^\beta \lambda(t)} u_{\leq A^\alpha \lambda(t)}$.\\
The estimate of the second term is given directly by the lemma.
\begin{eqnarray*}
 &&\|P_{> \lambda(t)A}(F(u_{\leq A^\beta \lambda(t)}))\|_{L^1 L^{\frac{6}{5-2s}}}\\
 &\lesssim& (\frac{\lambda(t) A^\beta}{\lambda(t)A})^{2-s} \|u_{\leq A^\beta \lambda(t)}\|_{L^p L^{3p}}^p \sup \|u_{\leq A^\beta \lambda(t)}\|_{\dot{H}^s}\\
 &\lesssim& A^{-(1-\beta)(2-s)}.
\end{eqnarray*} 
Combining these two estimates and taking sup for all time $t$, we can conclude the inequality (\ref{re1}).\\

\paragraph{Proof of the Second Inequality} To prove the inequality (\ref{re2}) we first define $t_i$ for $i \geq 1$ given $t_0 \in I$. 
\begin{equation}
 t_i = t_{i-1} + d \lambda^{-1} (t_{i-1}).
\label{boxdef}
\end{equation}
By the choice of $d$, all $t_i$'s are in the maximal lifespan $I$. See (\ref{choiced}) for more details. Please note 
that in the following argument we only need the cases $i=0,1,2$. But the definition of $t_i$ for all positive integers $i$ will be
used in later sections.\\
By the Strichartz estimate and the Duhamel formula, we have 
\begin{eqnarray*}
 &&\|u_{> \lambda(t_0) A}\|_{X([t_0,t_1])}\\
 &=& \left\|\int_t^{\infty} \frac{\sin((s-t)\sqrt{-\Delta})}{\sqrt{-\Delta}} P_{> \lambda(t_0)A} F(u(s)) ds\right\|_{X([t_0,t_1])}\\
 &\lesssim& \left\|\int_t^{t_2} \frac{\sin((s-t)\sqrt{-\Delta})}{\sqrt{-\Delta}} P_{> \lambda(t_0)A} F(u(s)) ds\right\|_{X([t_0,t_1])}\\
&& + \limsup_{T \rightarrow \infty} \left\|\int_{t_2}^{T} \frac{\sin((s-t)\sqrt{-\Delta})}{\sqrt{-\Delta}} P_{> \lambda(t_0)A} F(u(s)) ds\right\|_{X([t_0,t_1])}\\
 &\lesssim& \|P_{> \lambda(t_0) A} F(u(s))\|_{L^1 L^{\frac{6}{5-2s}}([t_0,t_2] \times \Rm^3)}\\
&& + \limsup_{T \rightarrow \infty} \left\|\int_{t_2}^{T} \frac{\sin((s-t)\sqrt{-\Delta})}{\sqrt{-\Delta}} P_{> \lambda(t_0)A} F(u(s))ds\right\|_{X([t_0,t_1])}\\
&=& I_1 + I_2. 
\end{eqnarray*}
The first term can be dominated by 
\begin{eqnarray*}
 I_1 &\lesssim& \|P_{> \lambda(t_0) A} F(u(s))\|_{L^1 L^{\frac{6}{5-2s}}([t_0,t_1] \times \Rm^3)}\\
&& + \|P_{> \lambda(t_0) A} F(u(s))\|_{L^1 L^{\frac{6}{5-2s}}([t_1,t_2] \times \Rm^3)}\\
&\lesssim& N(A) + N(\frac{\lambda(t_0)}{\lambda(t_1)}A)\\
&\lesssim& N(A^{1-\eps_1})
\end{eqnarray*}
for any small positive number $\eps_1$ and sufficiently large $A > A_0 (u,\eps_1)$, because $\lambda(t_0)$ and $\lambda(t_1)$ are comparable by the argument in the earlier sections.\\

For the second term, we will first find an upper bound of
\[
\left\|\int_{t_2}^T \frac{\sin((s-t)\sqrt{-\Delta})}{\sqrt{-\Delta}} F(u(s))ds\right\|_{L^\infty L^\infty}
\]
and then use an interpolation argument. If $x$ is small, we have 
\begin{eqnarray*}
 && \left|\left(\int_{t_2}^T \frac{\sin((s-t)\sqrt{-\Delta})}{\sqrt{-\Delta}} F(u(s))ds\right)(x)\right|\\
&=& \left|\int_{t_2}^T \int_{|y-x| = s-t} \frac{1}{4\pi (s-t)} F(u(s,y)) dS(y)ds\right|\\
&\lesssim& \int_{t_2}^T \int_{|y-x| = s-t} \frac{1}{4\pi (s-t)} |u(s,y)|^{p+1} dS(y)ds\\
&\lesssim& \int_{t_2}^T \int_{|y-x| = s-t} \frac{1}{(s-t)} \frac{1}{|y|^{(2/p)(p+1)}} dS(y)ds.
\end{eqnarray*}
In the last step, we use the following estimate for radial $\dot{H}^s$ functions. (Please see
lemma 3.2 of \cite{km}) 
\[
 |u(y)| \lesssim \frac{1}{|y|^{2/p}} \|u\|_{\dot{H}^s}.
\]
If $|x| \leq \frac{1}{2} (t_2 - t_1) \sim \lambda^{-1}(t_1) \sim \lambda^{-1}(t_0)$, then on the sphere for the integral
\[
 |y| \geq |s-t| - |x| \geq \frac{1}{2} (s-t).
\]
Thus for these small $x$, 
\begin{eqnarray*}
 && \left|\left(\int_{t_2}^T \frac{\sin((s-t)\sqrt{-\Delta})}{\sqrt{-\Delta}} F(u(s))ds\right)(x)\right|\\
&\lesssim& \int_{t_2}^T \int_{|y-x| = s-t} \frac{1}{(s-t)} \frac{1}{(s-t)^{(2/p)(p+1)}} dS(y)ds\\
&\lesssim& \int_{t_2}^T \int_{|y-x| = s-t} \frac{1}{(s-t)^{3 + 2/p}} dS(y)ds\\
&\lesssim& \int_{t_2}^T \frac{(s-t)^2}{(s-t)^{3 + 2/p}} ds\\
&\lesssim& \int_{t_2}^T \frac{1}{(s-t)^{1 + 2/p}} ds\\
&\lesssim& (t_2 - t)^{-2/p}\\
&\lesssim& (t_2 - t_1)^{-2/p} \sim [\lambda(t_0)]^{2/p}.
\end{eqnarray*}
On the other hand, we also have a uniform bound for all $t \leq T' \leq T$ using our assumption (\ref{add1}).
\begin{equation}
 \left\|\left(\begin{array}{l} \displaystyle \int_{T'}^T \frac{\sin((s-t)\sqrt{-\Delta})}{\sqrt{-\Delta}} F(u(s)) ds\\
\displaystyle -\int_{T'}^T \cos((s-t)\sqrt{-\Delta}) F(u(s)) ds \end{array} \right)\right\|_{\dot{H}^s \times \dot{H}^{s-1}} \lesssim_u 1.
\label {uniformhs}
\end{equation} 
This gives us an estimate for large $x$. If $|x| > \frac{1}{2} (t_2 - t_1)$, we have  
\begin{eqnarray*}
 && \left|\left(\int_{t_2}^T \frac{\sin((s-t)\sqrt{-\Delta})}{\sqrt{-\Delta}} F(u(s))ds\right)(x)\right|\\
 &\lesssim& \frac{1}{|x|^{2/p}} \left\|\int_{t_2}^T \frac{\sin((s-t)\sqrt{-\Delta})}{\sqrt{-\Delta}} F(u(s)) ds\right\|_{\dot{H}^s}\\
&\lesssim& \frac{1}{(t_2 - t_1)^{2/p}} \\
&\simeq& [\lambda(t_0)]^{2/p}.
\end{eqnarray*}
Combining the estimates for small and large $x$, we obtain
\begin{equation}
 \left\|\int_{t_2}^T \frac{\sin((s-t)\sqrt{-\Delta})}{\sqrt{-\Delta}} F(u(s))ds\right\|_{L^\infty L^\infty} \lesssim \lambda(t_0)^{2/p}.
 \label{linftybound}
\end{equation}
This implies 
\begin{equation}
 \left\|P_{> \lambda(t_0)A} \int_{t_2}^T \frac{\sin((s-t)\sqrt{-\Delta})}{\sqrt{-\Delta}} F(u(s))ds\right\|_{L^\infty L^\infty} \lesssim \lambda(t_0)^{2/p}.
\end{equation}
By (\ref{uniformhs}), we also have 
\begin{equation}
 \left\|P_{> \lambda(t_0)A} \int_{t_2}^T \frac{\sin((s-t)\sqrt{-\Delta})}{\sqrt{-\Delta}} F(u(s))ds\right\|_{L^\infty L^2} \lesssim (\lambda(t_0)A)^{-s}.
\end{equation}
Thus 
\begin{eqnarray*}
 && \left\|P_{> \lambda(t_0)A} \int_{t_2}^T \frac{\sin((s-t)\sqrt{-\Delta})}{\sqrt{-\Delta}} F(u(s))ds\right\|_{L^\infty L^{{3p}/2}}\\
&\leq& \|\cdotp\|_{L^\infty L^\infty}^{\displaystyle 1 - 4/{3p}} \|\cdotp\|_{L^\infty L^2}^{\displaystyle 4/{3p}}\\
&\lesssim& [\lambda(t_0)^{\displaystyle 2/p}]^{\displaystyle 1 - 4/{3p}} [(\lambda(t_0)A)^{-s}]^{\displaystyle 4/{3p}}\\
&=& A^{\displaystyle -{4s}/{3p}}.
\end{eqnarray*}
Thus 
\begin{eqnarray*}
 && \left\|P_{> \lambda(t_0)A} \int_{t_2}^T \frac{\sin((s-t)\sqrt{-\Delta})}{\sqrt{-\Delta}} F(u(s))ds\right\|_{X([t_0,t_1])}\\
 &\lesssim& \|\cdotp\|_{\displaystyle L^\infty L^{{3p}/2}([t_0,t_1]\times \Rm^3)}^{\displaystyle 1 - {2}/{s(p+1)}} \|\cdotp\|_{\displaystyle L^{2/s} L^{2/{1-s}}([t_0,t_1]\times \Rm^3)}^{\displaystyle {2}/{s(p+1)}}\\
 &\lesssim& A^{\displaystyle ({-4s}/{3p}) (1 -{2}/{s(p+1)})} \\
 &\lesssim& A^{\displaystyle - \frac{4}{3p} (s - \frac{2}{p+1})}\\
 &\lesssim& A^{\displaystyle - \frac{4(sp + s -2)}{3p(p+1)}}. 
\end{eqnarray*}
In the second step we use the fact that the $L^{2/s} L^{2/{1-s}}$ norm is uniformly bounded. This comes from the uniform 
$\dot{H}^s \times \dot{H}^{1-s}$ bound at $t = t_1$ in (\ref{uniformhs}) and the Strichartz estimate. 
(Note that $(2/s,2/{1-s})$ is an $s$-admissible pair)\\

Letting $T \rightarrow \infty$, we have the estimate for $I_2$. This completes the proof of (\ref{re2}).
\paragraph{The Decay of $N(A)$ and $S(A)$} Let us assume $s >9/10$. Then by (\ref{sigma}), $\sigma$ in the recurrence inequality can be any positive real number less than
\[
 3 \min \left\{\frac{1}{2} - \frac{1}{p+1} - \frac{5-2s}{6(p+1)}, \frac{5 -2s}{6(p+1)}\right\}.
\]
The second bound is a decreasing function in $s$. So this bound is greater than $3/10$, which is the value of the function when $s=1$ and $p=4$.
The first bound is an increasing function of $s$. Thus it is greater than the value when $s=9/10$. If $s=9/10$, we have $p=10/3$. Thus
\[
 \frac{1}{2} - \frac{1}{p+1} - \frac{5-2s}{6(p+1)} = \frac{1}{2} - \frac{1}{10/3 + 1} - \frac{5 -9/5}{6(10/3 +1)} = 19/130 > 1/10.
\]
Thus we can choose $\sigma = 3/10$.\\
Now let us look at the value of $\sigma_1$. If $s>9/10$, then the numerator of $\sigma_1$ is greater than
\[
 4(\frac{9}{10} \frac{10}{3} + \frac{9}{10} -2) = 76/10.
\]
While its denominator $3p(p+1)$ is less than $3 \times 4 \times 5 = 60$. Thus $\sigma_1 > 1/8$. This gives us the recurrence inequalities for $s >9/10$
\begin{equation}
 N(A) \lesssim S(A^\beta) S^{p}(A^\alpha) + A^{-\frac{3}{10}(\beta - \alpha)} + A^{-(2-s)(1-\beta)}.
 \label{re3}
\end{equation}
\[
 S(A) \lesssim N(A^{1-\eps_1}) + A^{-1/8}.
\]
For each sufficiently large $A$, plug the first inequality into the second one, we have 
\begin{eqnarray*}
 S(A) \lesssim S(A^{(1-\eps_1)\beta}) S^{p}(A^{(1-\eps_1)\alpha}) &+&  A^{-\frac{3}{10}(1-\eps_1)(\beta - \alpha)}\\
 &+& A^{-(2-s)(1-\eps_1)(1-\beta)} + A^{-1/8}.
\end{eqnarray*}
Choose $\alpha$, $\beta$ and $\eps_1$ so that 
\begin{equation}
 (1-\eps_1)\beta = 0.85;\,  (1-\eps_1)\alpha = 0.4;\,  \eps_1=1/10000.
\end{equation}
Then we have 
\[
 S(A) \lesssim S(A^{0.85})S^p (A^{0.4}) + A^{-1/8}
\]
with the additional information that $S(A) \rightarrow 0$ as $A \rightarrow \infty$. Using the following lemma, we have $S(A) \lesssim A^{-1/8}$.
\paragraph{Lemma} Suppose $S(A) \rightarrow 0$ as $A \rightarrow \infty$. In addition, there exist $\alpha,\beta \in (0,1)$ and $p, \omega > 0$ with
\[
 p \alpha + \beta > 1,
\]
such that  
\[
 S(A) \lesssim S(A^\beta)S^p (A^\alpha) + A^{-\omega}
\]
is true for each sufficiently large $A$. Then 
\[
 S(A) \lesssim A^{-\omega}.
\]
for each sufficiently large $A$. \\
\newline 
Now we come back to (\ref{re3}) with $S(A) \lesssim A^{-1/8}$. By our choice of $\alpha$ and $\beta$, we have 
\[
 N(A) \lesssim A^{-0.13}
\]
for each sufficiently large $A$. Observing that our $N(A)$ and $S(A)$ are uniformly bounded, we know that the decay inequalities above are true for all $A$.
\paragraph{Short-time Contribution in the Duhamel Formula} Define 
\begin{equation}
 v_{t'}({\tilde{t}}) = \int_{t'}^{t' + d\lambda^{-1}(t')} \frac{\sin((s-\tilde{t})\sqrt{-\Delta})}{\sqrt{-\Delta}} F(u(s)) ds
\label{defs}
\end{equation}
Thus
\[
 \partial_t v_{t'}(\tilde{t}) = -\int_{t'}^{t' + d\lambda^{-1}(t')} \cos((s-\tilde{t})\sqrt{-\Delta}) F(u(s)) ds.
\]
for all $\tilde{t}\leq t'$. This is a short time contribution in the Duhamel formula. One can also think it to be the solution of the backward time problem
\[
 \partial_{tt} v - \Delta v = \chi([t', t'+ d\lambda^{-1}(t')]) F(u(s))
\]
with the initial data $(0,0)$ at time $t' + d \lambda^{-1}(t')$. Here $\chi$ is the characteristic function of the time period indicated.
By the Strichartz estimate, we have 
\begin{eqnarray*}
 && \|P_{> \lambda(t')A} (v_{t'}(\tilde{t}), \partial_t v_{t'}(\tilde{t}))\|_{\dot{H}^s \times \dot{H}^{s-1}}\\
 &\lesssim& \|P_{> \lambda(t')A}(F(u))\|_{\displaystyle L^1 L^{\frac{6}{5-2s}}([t', t'+ d\lambda^{-1}(t')]\times \Rm^3)}\\
 &\lesssim& N(A) \lesssim A^{-0.13}.
\end{eqnarray*}
This implies 
\begin{eqnarray*}
 && \|P_{ \lambda(t')A < \cdotp < 2\lambda(t')A} (v_{t'}(\tilde{t}), \partial_t v_{t'}(\tilde{t}))\|_{\dot{H}^{s +1/8} \times \dot{H}^{s +(1/8)-1}}\\
 &\lesssim& (\lambda(t')A)^{1/8} \|P_{ \lambda(t')A < \cdotp < 2\lambda(t')A} (v_{t'}(\tilde{t}), \partial_t v_{t'}(\tilde{t}))\|_{\dot{H}^{s} \times \dot{H}^{s -1}}\\
&\lesssim& (\lambda(t')A)^{1/8} \|P_{> \lambda(t')A} (v_{t'}(\tilde{t}), \partial_t v_{t'}(\tilde{t}))\|_{\dot{H}^s \times \dot{H}^{s-1}}\\
&\lesssim& (\lambda(t')A)^{1/8} A^{-0.13}\\
&\lesssim& \lambda(t')^{1/8} A^{-1/200}.
\end{eqnarray*}
Take $A = 2^k$ for $k \geq 0$ and sum. We obtain
\[
 \|P_{ > \lambda(t')} (v_{t'}(\tilde{t}), \partial_t v_{t'}(\tilde{t}))\|_{\dot{H}^{s +1/8} \times \dot{H}^{s +(1/8)-1}} \lesssim \lambda(t')^{1/8}.
\]
We also have the low frequency estimate 
\begin{eqnarray*}
 && \|P_{ \leq \lambda(t')} (v_{t'}(\tilde{t}), \partial_t v_{t'}(\tilde{t}))\|_{\dot{H}^{s +1/8} \times \dot{H}^{s +(1/8)-1}}\\
 &\lesssim& \lambda(t')^{1/8} 
\| (v_{t'}(\tilde{t}), \partial_t v_{t'}(\tilde{t}))\|_{\dot{H}^{s} \times \dot{H}^{s-1}} \lesssim \lambda(t')^{1/8}.
\end{eqnarray*}
Combining the high and low frequency parts we have 
\begin{equation}
 \|(v_{t'}(\tilde{t}), \partial_t v_{t'}(\tilde{t}))\|_{\dot{H}^{s +1/8} \times \dot{H}^{s +(1/8)-1}} \lesssim \lambda(t')^{1/8}.
\label{localreg}
\end{equation}
for all $\tilde{t} \leq t'$. 
\subsection{Additional Regularity} In this section, we will show additional regularity of the solution, at least locally in space. 
The idea is to estimate the local $\dot{H}^1 \times L^2$ norm of $(u(t_0), \partial_t u(t_0))$ by separating the time interval $[t_0, +\infty)$ into\\
(i) finitely many intervals corresponding to the boxes below (using the short time estimate we obtained above), plus\\
(ii) infinitely many intervals corresponding to the thin slices below (long time estimate),\\
using the Duhamel formula. The strong Huygens' principle plays an important role in this argument. Let us first construct the boxes and slices mentioned above. 
\paragraph{The Construction of Boxes} Let us fix $t_0 = 0 \in I$. Define $t_i$ as before
\[
 t_i = t_{i-1} + d \lambda^{-1} (t_{i-1}).
\]
The $i$'th box is the circular cylinder in the space-time
\[
 B_i = \{(x,t): |x| \leq R_0 \lambda^{-1}(t_i), t \in [t_i,t_{i+1}] \}.
\]
Here the constant $R_0 = R_0(u)$ is the same constant that appeared in the local compactness section, part (iv).
For $i \geq 1$, we have 
\[
 \lambda^{-1}(t_i) \leq C_1 \lambda^{-1}(t_{i-1}) = C_1 \frac{t_i - t_{i-1}}{d} \leq \frac{C_1}{d} (t_i-t_0).
\]
Thus for all $(x,t) \in B_i$ with $i \geq 1$, we have 
\begin{equation}
 |x| + t - t_0 \leq (\frac{R_0 C_1}{d} + C_1 + 1) (t_i -t_0).
\label{boxestimate}
\end{equation}
All the constants mentioned here are from the local compactness part, thus they do not depend on the integer $i$. This inequality will be useful later. 
Given a time $T > 0$, let us define 
\[
 m(T) = \sup\{m \in {\mathbb Z}: t_m \in [t_0,t_0 + T]\}.
\]
Thus $m(T)$ is a nonnegative integer. $m(T)$ can never be infinity because the frequency
scale function $\lambda(t)$ is bounded.(Please see the minimal blow-up solution part)
This function helps us determine how many boxes to use before we turn to slices. 
\paragraph{Construction of Slices} The slices begin at $t_{m(T)+1}$, where the last box ends. Let us call 
\[
 T_0 (T) = t_{m(T)+1}.
\]
and define 
\[
 T_{i+1} (T) = \eta (T_{i} (T) - t_0) + T_i (T) = (1+\eta)^{i+1} (T_0 (T)-t_0) + t_0.
\]
Here the small constant $\eta$ should be chosen so that 
\begin{equation}
 0< \eta < \min \left\{\frac{d}{2 R_1}, 1/100\right\}. 
\label{choiceofeta}
\end{equation}
The constants in the definition come from the local compactness part, so they depend only on $u$.    
\paragraph{Local Estimate of $u(t_0)$} Now we will estimate $u(t_0)$ locally in space. 
Namely, we will do the estimate in a ball with radius $R$. Let us choose $R > d \lambda^{-1}(t_0)$ and $T = 100R$.
By the Duhamel formula, 
\begin{eqnarray*}
 u(t_0) &=& v_{t_0} (t_0) + v_{t_1} (t_0) + v_{t_2} (t_0) + \cdots + v_{t_{m(100R)}}(t_0)\\
 && + \bar{v}_0 + \bar{v}_1 + \bar{v}_2 + \cdots,\\
 \partial_t  u(t_0) &=& \partial_t v_{t_0} (t_0) + \partial_t v_{t_1} (t_0) + \partial_t v_{t_2} (t_0) + \cdots + \partial_t v_{t_{m(100R)}}(t_0)\\
 && + \partial_t \bar{v}_0 + \partial_t \bar{v}_1 + \partial_t \bar{v}_2 + \cdots 
 \end{eqnarray*}
as a weak limit in $\dot{H}^s \times \dot{H}^{s-1}$. Here the $v_{t_i}(t_0)$'s are defined 
in (\ref{defs}) while the $\bar{v_i}$'s are the contribution from the thin slices.
\begin{equation}
\bar{v}_i = \int_{T_i}^{T_{i+1}} \frac{\sin((s-t_0)\sqrt{-\Delta})}{\sqrt{-\Delta}} F(u(s)) ds
 \label{defbar}
\end{equation}
\[
 \partial_t \bar{v}_i = -\int_{T_i}^{T_{i+1}} \cos((s-t_0)\sqrt{-\Delta}) F(u(s)) ds.
\]
In order to use the decay of $u(x, t)$ as $x$ is large, we also define 
\begin{equation}
\tilde{v}_i = \int_{T_i}^{T_{i+1}} \frac{\sin((s-t_0)\sqrt{-\Delta})}{\sqrt{-\Delta}} [\chi_{|x|> {3(T_i-t_0)}/4} F(u(s))] ds
 \label{deftilde}
\end{equation}
\[
 \partial_t \tilde{v}_i = -\int_{T_i}^{T_{i+1}} \cos((s-t_0)\sqrt{-\Delta}) [\chi_{|x|>{3(T_i-t_0)}/4} F(u(s))] ds. 
\]
Here the function $\chi$ is a cutoff function so that we discard the center part of the nonlinearity. Because the functions $u(t_0)$ and $\partial_t u(t_0)$ in the ball of radius $R$ only depend on the nonlinearity in the region
\[
 \{(x,t): t-t_0-R \leq |x| \leq t-t_0 +R\} 
\]
in the Duhamel formula, we know this part of $(u(t_0),\partial_t u(t_0))$ is not affected by the cutoff. Thus we have 
\begin{eqnarray}
 u(t_0) & = & v_{t_0} (t_0) + v_{t_1} (t_0) + v_{t_2} (t_0) + \cdots + v_{t_{m(100R)}}\nonumber\\
 & & + \tilde{v}_0 + \tilde{v}_1 + \tilde{v}_2 + \cdots. \label{sum}\\
 \partial_t  u(t_0) &=& \partial_t v_{t_0} (t_0) + \partial_t v_{t_1} (t_0) + \partial_t v_{t_2} (t_0) + \cdots + \partial_t v_{t_{m(100R)}}\nonumber\\
 && + \partial_t \tilde{v}_0 + \partial_t \tilde{v}_1 + \partial_t \tilde{v}_2 + \cdots\label{sum1}
\end{eqnarray}
as a weak limit in the ball $B(0,R)$. 
\paragraph{Short-time Contribution} By (\ref{localreg}), we have 
\begin{eqnarray*}
 &&\|v_{t_0} (t_0) + v_{t_1} (t_0) + v_{t_2} (t_0) + \cdots + v_{t_{m(100R)}}(t_0)\|_{\dot{H}^{s+1/8}} \\
&& + \|\partial_t v_{t_0} (t_0) + \partial_t v_{t_1} (t_0) + \partial_t v_{t_2} (t_0) + \cdots + \partial_t v_{t_{m(100R)}}(t_0)\|_{\dot{H}^{s +1/8 -1}}\\
 &\lesssim& \lambda(t_0)^{1/8} + \lambda(t_1)^{1/8} + \cdots + \lambda(t_{m(100R)})^{1/8}.
\end{eqnarray*}
Let us define 
\begin{equation}
 Q(T) = \lambda(t_0)^{1/8} + \lambda(t_1)^{1/8} + \cdots + \lambda(t_{m(T)})^{1/8}.
 \label{defQ}
\end{equation}
Combining this estimate with the uniform $\dot{H}^s \times \dot{H}^{s-1}$ bound of the short time contribution we have 
\begin{eqnarray*}
 &&\|v_{t_0} (t_0) + v_{t_1} (t_0) + v_{t_2} (t_0) + \cdots + v_{t_{m(100R)}}\|_{\dot{H}^{1}} \\
&& + \|\partial_t v_{t_0} (t_0) + \partial_t v_{t_1} (t_0) + \partial_t v_{t_2} (t_0) + \cdots + \partial_t v_{t_{m(100R)}}\|_{L^2}\\ 
&\lesssim& Q(100R)^{8(1-s)}.
\end{eqnarray*}
and 
\begin{eqnarray*}
 &&\|v_{t_0} (t_0) + v_{t_1} (t_0) + v_{t_2} (t_0) + \cdots + v_{t_{m(100R)}}(t_0)\|_{\displaystyle \dot{H}^{\frac{3p}{2(p+2)}}} \\
&& + \|\partial_t v_{t_0} (t_0) + \partial_t v_{t_1} (t_0) + \partial_t v_{t_2} (t_0) + \cdots + \partial_t v_{t_{m(100R)}}(t_0)\|_{\displaystyle \dot{H}^{\frac{3p}{2(p+2)}-1}}\\ 
&\lesssim& Q(100R)^{\displaystyle 8(\frac{3p}{2(p+2)}-s)}\\
&=& Q(100R) ^ { \displaystyle \frac{16(1-s)}{p+2}}.
\end{eqnarray*}
By the Sobolev embedding, this implies
\begin{equation}
 \|v_{t_0} (t_0) + v_{t_1} (t_0) + v_{t_2} (t_0) + \cdots + v_{t_{m(100R)}}(t_0)\|_{L^{p+2}} \lesssim Q(100R) ^ {\frac{16(1-s)}{p+2}}.
 \label{shorttime1}
\end{equation}

\paragraph{Long-time Contribution} By the Strichartz estimate
\begin{eqnarray*}
 &&\|(\tilde{v}_i, \partial_t \tilde{v}_i)\|_{\dot{H}^1 \times L^2}\\
 &\lesssim& \|\chi_{|x|>{3(T_i-t_0)}/4} F(u(s))\|_{L^1 L^2 ([T_i, T_{i+1}] \times \Rm^3)}\\
 &\lesssim& \|\chi u(s)\|_{L^p L^{3p}([T_i, T_{i+1}] \times \Rm^3)}^p \|\chi u(s)\|_{L^\infty L^6 ([T_i, T_{i+1}] \times \Rm^3)}.
\end{eqnarray*}
Similarly 
\begin{eqnarray*}
 &&\|(\tilde{v}_i, \partial_t \tilde{v}_i)\|_{\dot{H}^{\frac{3p}{2(p+2)}} \times \dot{H}^{\frac{3p}{2(p+2)}-1}}\\
 &\lesssim& \|\chi_{|x|>{3(T_i-t_0)}/4} F(u(s))\|_{\displaystyle L^1 L^{\frac{3(p+2)}{p+5}} ([T_i, T_{i+1}] \times \Rm^3)}\\
 &\lesssim& \|\chi u(s)\|_{L^p L^{3p}([T_i, T_{i+1}] \times \Rm^3)}^p \|\chi u(s)\|_{L^\infty L^{p+2} ([T_i, T_{i+1}] \times \Rm^3)}.
\end{eqnarray*}
We will show that the first factor is uniformly bounded for all $i$, while the second factor in each estimate decays so that we can take a sum of all $\tilde{v}_i$'s.
\paragraph{Boundedness of the First Factor} There are two cases.\\
(I) Case 1. If $\lambda(T_i) \geq 2 R_1 (T_i - t_0)^{-1}$, then by local compactness part (iii), we have 
\[
 \left\|\left(\varphi(\frac{x}{R_1 \lambda^{-1}(T_i)})u(T_i),\varphi(\frac{x}{R_1 \lambda^{-1}(T_i)})\partial_t u(T_i)\right)\right\|_{\dot{H}^1 \times \dot{H}^{s-1}} \leq \delta.  
\]
Thus the solution $\bar{u}$ of the NLW with the initial data 
\[
\left(\varphi(\frac{x}{R_1 \lambda^{-1}(T_i)})u(T_i),\varphi(\frac{x}{R_1 \lambda^{-1}(T_i)})\partial_t u(T_i)\right)
\]
scatters and 
\[
 \|\bar{u}\|_{L^p L^{3p} (\Rm \times \Rm^3)} \lesssim 1. 
\]
By the finite speed of propagation and the definition of the cutoff function $\varphi$, we have 
\[
 \|u\|_{L^p L^{3p} (\Omega)} \lesssim 1.
\]
Here $\Omega$ is given by
\[
 \Omega = \{(x,t): |x| > R_1 \lambda^{-1}(T_i) + |t - T_i|\}.
\]
Let us check in case I, that the region
\[
 \{(x,t): |x| > 3(T_i-t_0)/4, t \in [T_i, T_{i+1}]\}
\]
is completely contained in $\Omega$. 
By the assumption, we have $R_1 \lambda^{-1}(T_i) \leq (1/2) (T_i - t_0)$. Thus if $|x| > (3/4) (T_i - t_0)$ and $t \in [T_i,T_{i+1}]$, we have 
\[
 |t - T_i| + R_1 \lambda^{-1}(T_i) \leq (T_{i+1}-T_i) + \frac{1}{2} (T_i - t_0) \leq (\eta + \frac{1}{2})(T_i-t_0) < |x|.
\]
Thus in case I
\[
 \|\chi u(s)\|_{L^p L^{3p}([T_i, T_{i+1}] \times \Rm^3)} \leq \|u\|_{L^p L^{3p} (\Omega)}  \lesssim 1. 
\]
(II) Case 2. If $\lambda(T_i) < 2 R_1 (T_i - t_0)^{-1}$, then
\[
 \lambda^{-1}(T_i) > \frac{T_i - t_0}{2 R_1}.
\]
This implies 
\[
 d \lambda^{-1} (T_i) > \frac{d}{2 R_1} (T_i - t_0) > \eta (T_i -t_0) = T_{i+1} - T_i.
\]
by the choice of $\eta$. So
\begin{equation}
 \|\chi u(s)\|_{L^p L^{3p}([T_i, T_{i+1}] \times \Rm^3)}
 \lesssim \|u(s)\|_{L^p L^{3p} ([T_i, T_i + d \lambda^{-1}(T_i)]\times \Rm^3)} \lesssim 1
\end{equation}
by the local compactness estimate part (ii).\\
In summary, the first factor is always uniformly bounded. 
\paragraph{Decay of the Second Factor} The estimate is straight forward, for all $t \in [T_i, T_{i+1}]$
\begin{eqnarray*}
 \int_{|x|>\frac{3}{4}(T_i - t_0)} |u(t)|^6 dx &\leq& \int_{|x|>\frac{3}{4}(T_i - t_0)} |u(t)|^{3p/2} dx \\
& & \times \left(\sup_{|x|>\frac{3}{4}(T_i - t_0)}|u(t)|\right)^{6 - {3p}/2}\\
&\lesssim& \int |u(t)|^{3p/2} dx \left(\frac{\|u(t)\|_{\dot{H}^s}}{[\frac{3}{4} (T_i-t_0)]^{2/p}}\right)^{6 - 3p/2}\\
&\lesssim& \frac{1}{(T_i - t_0)^{12/p - 3}}\\
&\lesssim& \frac{1}{(T_0 - t_0)^{12/p - 3}} \left(\frac{1}{(1 + \eta)^{12/p -3}}\right)^i.
\end{eqnarray*}
Thus 
\[
 \|\chi u(s)\|_{L^\infty L^6 ([T_i, T_{i+1}] \times \Rm^3)} \lesssim \frac{1}{(T_0 -t_0)^{2/p -1/2}} \left(\frac{1}{(1 + \eta)^{2/p -1/2}}\right)^i.
\]
We know $2/p - 1/2 = 1 -s $. So
\[
 \|\chi u(s)\|_{L^\infty L^6 ([T_i, T_{i+1}] \times \Rm^3)} \lesssim \frac{1}{(T_0 -t_0)^{1-s}} \left(\frac{1}{(1 + \eta)^{1-s}}\right)^i.
\]
Similarly
\begin{eqnarray*}
\int_{|x|>\frac{3}{4}(T_i - t_0)} |u(t)|^{p+2} dx &\leq& \int_{|x|>(3/4)(T_i - t_0)} |u(t)|^{3p/2} dx \\
&& \times \left(\sup_{|x|>\frac{3}{4}(T_i - t_0)}|u(t)|\right)^{(p+2) - {3p}/2}\\
&\lesssim& \int |u(t)|^{3p/2} dx \left(\frac{\|u(t)\|_{\dot{H}^s}}{[\frac{3}{4} (T_i-t_0)]^{2/p}}\right)^{2 - p/2}\\
&\lesssim& \frac{1}{(T_i - t_0)^{4/p - 1}}\\
&=& \frac{1}{(T_i - t_0)^{2(1-s)}}\\
&\lesssim& \frac{1}{(T_0 - t_0)^{2(1-s)}} \left(\frac{1}{(1 + \eta)^{2(1-s)}}\right)^i.
\end{eqnarray*}
Thus
\[
 \|\chi u(s)\|_{L^\infty L^{p+2} ([T_i, T_{i+1}] \times \Rm^3)} \lesssim \frac{1}{(T_0 -t_0)^{\frac{2(1-s)}{p+2}}} \left(\frac{1}{(1 + \eta)^{\frac{2(1-s)}{p+2}}}\right)^i.
\]

\paragraph{The End of the Long-time Contribution} Combining the estimates for the two factors, we have 
\begin{equation}
 \|(\tilde{v}_i, \partial_t \tilde{v}_i)\|_{\dot{H}^1 \times L^2} \lesssim \frac{1}{(T_0 -t_0)^{1-s}} \left(\frac{1}{(1 + \eta)^{1-s}}\right)^i.
 \label{longtimereg}
\end{equation}
and 
\[
 \|(\tilde{v}_i, \partial_t \tilde{v}_i)\|_{\dot{H}^{\frac{3p}{2(p+2)}} \times \dot{H}^{\frac{3p}{2(p+2)}-1}} \lesssim \frac{1}{(T_0 -t_0)^{\frac{2(1-s)}{p+2}}} \left(\frac{1}{(1 + \eta)^{\frac{2(1-s)}{p+2}}}\right)^i.
\]
By the Sobolev embedding the second estimate implies
\begin{equation}
 \|\tilde{v}_i\|_{L^{p+2}} \lesssim \frac{1}{(T_0 -t_0)^{\frac{2(1-s)}{p+2}}} \left(\frac{1}{(1 + \eta)^{\frac{2(1-s)}{p+2}}}\right)^i.
 \label{longtime1}
\end{equation}
The estimate (\ref{longtimereg}) means that the pair consisting of the right hands of (\ref{sum}) and (\ref{sum1}) converges to some pair $(\tilde{u}_0, \tilde{u}_1)$ in $\dot{H}^1 \times L^2$ with the following estimate.
\begin{equation}
 \|(\tilde{u}_0, \tilde{u}_1)\|_{\dot{H}^1 \times L^2} \lesssim Q(100R)^{8(1-s)} + \frac{1}{(T_0 -t_0)^{1-s}}.
\end{equation}
By (\ref{longtime1}) and (\ref{shorttime1}), we also have 
\[
  \|\tilde{u}_0\|_{L^{p+2}} \lesssim Q(100R) ^ {\frac{16(1-s)}{p+2}} + \frac{1}{(T_0 -t_0)^{\frac{2(1-s)}{p+2}}}.
\]
We have $T_0 - t_0 = t_{m(100R) +1} - t_0 > 100R + t_0 - t_0 = 100 R$ by the definition of $m(T)$. Thus
\[
 \|(\tilde{u}_0, \tilde{u}_1)\|_{\dot{H}^1 \times L^2} \lesssim Q(100R)^{8(1-s)} + \frac{1}{(100 R)^{1-s}}.
\]
\[
 \|\tilde{u}_0\|_{L^{p+2}} \lesssim Q(100R) ^ { \frac{16(1-s)}{p+2}} + \frac{1}{(100R)^{\frac{2(1-s)}{p+2}}}.
\]
Because $R > d \lambda^{-1}(t_0)$, we have 
\[
 \frac{1}{100 R} \lesssim \lambda(t_0)  = (\lambda(t_0)^{1/8})^8 \leq Q(100R)^8.
\]
So we have 
\begin{equation}
 \|(\tilde{u}_0, \tilde{u}_1)\|_{\dot{H}^1 \times L^2} \lesssim Q(100R)^{8(1-s)}.
\label{estimatetilde1}
\end{equation}
\begin{equation}
 \|\tilde{u}_0\|_{L^{p+2}} \lesssim Q(100R) ^ { \frac{16(1-s)}{p+2}}.
\label{estimatetilde2}
\end{equation}
\paragraph{The Identity in the Ball} 
Now, in the ball of radius $R$, the pair that consists of the right hands of (\ref{sum}) and (\ref{sum1}) converges to 
$(\tilde{u}_0, \tilde{u}_1)$ strongly in $\dot{H}^1 \times L^2$, 
and to $(u(t_0), \partial_t u(t_0))$ weakly. Thus in the ball of radius $R$, 
\[
 u(t_0) = \tilde{u}_0; \,\, \partial_t u(t_0) = \tilde{u}_1.
\]
\section{The Death of Solutions}
We will show a contradiction in this section as $s$ is sufficiently close to $1$. 
There are two different cases.\\
(I) The function $Q(T)$ is bounded.\\
(II) The function $Q(T) \rightarrow \infty$ as $T \rightarrow \infty$.\\
The first case gives us a $\dot{H}^1 \times L^2$ estimate in the 
whole space, thus it is much easier to deal with. In either case, we will use the global
integral estimate
\[
  \int_0^T \int_{\Rm^3} \frac{|\tilde{u}|^{p+2}}{|x|} dx dt \leq \frac{p+2}{p} E.
\]
This is the only place where we use the defocusing condition. All arguments before this point are 
also valid in the focusing case. 
\subsection{Case 1} Let us assume $Q(T)$ is bounded. In other words,
\begin{equation}
 \lambda(t_0)^{1/8} + \lambda(t_1)^{1/8} + \lambda(t_2)^{1/8} + \cdots  \leq C. 
\label{boundofQ}
\end{equation}
By the Duhamel formula, we have 
\begin{eqnarray*}
 u(t_n) &=& v_{t_n} (t_n) + v_{t_{n+1}} (t_n) + v_{t_{n+2}} (t_n) + \cdots.\\
 \partial_t u(t_n) &=& \partial_t v_{t_n} (t_n) + \partial_t v_{t_{n+1}} (t_n) + \partial_t v_{t_{n+2}} (t_n) + \cdots.
\end{eqnarray*}
as a weak limit. At the same time the right hand has a strong limit in $\dot{H}^{s+1/8} \times \dot{H}^{(s+1/8) -1}$ 
by the estimate (\ref{localreg}) and the assumption (\ref{boundofQ}). Thus we have 
$(u(t_n), \partial_t u(t_n)) \in \dot{H}^{s+1/8} \times \dot{H}^{(s+1/8) -1}$ with norm
\[
 \|(u(t_n), \partial_t u(t_n))\|_{\dot{H}^{s+1/8} \times \dot{H}^{(s+1/8) -1}} \lesssim \lambda(t_n)^{1/8} + \lambda(t_{n+1})^{1/8} + \cdots.
\]
This implies 
\[
 \lim_{n \rightarrow \infty} \|(u(t_n), \partial_t u(t_n))\|_{\dot{H}^{s+1/8} \times \dot{H}^{(s+1/8) -1}} = 0.
\]
Using the interpolation between $\dot{H}^{s +1/8}$ and $\dot{H}^s$, we have 
\[
 \lim_{n \rightarrow \infty} \|(u(t_n), \partial_t u(t_n))\|_{\dot{H}^{1} \times L^2} = 0.
\]
Thus by the Sobolev embedding, 
\[
 \lim_{n \rightarrow \infty} \|u(t_n)\|_{L^6} = 0.
\]
Since $u(t)$ is uniformly bounded in $\dot{H}^s$, it is also uniformly bounded in $L^{3p/2}$ by the Sobolev embedding. 
By the inequality $3p/2 < p+2 < 6$, we have 
\[
 \lim_{n \rightarrow \infty} \|u(t_n)\|_{L^{p+2}} = 0.
\]
Thus the energy at time $t_n$ converges to zero. 
\[
 E(t_n) = \int_{\Rm^3} (|\nabla_x u|^2 + |\partial_t u|^2 + 2 \frac{|u(x)|^{p+2}}{p+2}) dx \rightarrow 0.
\]
This implies (Since the wave equation is time reversible)
\[
 \int_{t_0}^{t_n} \int_{\Rm^3} \frac{|u(x,s)|^{p+2}}{|x|} dx ds \lesssim E(t_n) \rightarrow 0. 
\]
Letting $n \rightarrow \infty$, we have $u \equiv 0$. This is a contradiction. 
\subsection{Case 2}
In this case let us assume $Q(T) \rightarrow \infty$. We need to use our 
local estimate obtained in the previous section. By (\ref{estimatetilde1}) and (\ref{estimatetilde2})
\[
 E(\tilde{u}_0, \tilde{u}_1) \lesssim Q(100R)^{16(1-s)}.
\]
Let $\tilde{u}$ be the solution of the nonlinear wave equation with the initial data $(\tilde{u}_0, \tilde{u}_1)$.
In the defocusing, energy subcritical case with a finite energy, the solution will never break down in finite time, thus
\[
  \int_0^\infty \int_{\Rm^3} \frac{|\tilde{u}|^{p+2}}{|x|} dx dt \leq \frac{p+2}{p} E \lesssim Q(100R)^{16(1-s)}.
\]
Observing in the ball of radius $R$ this pair of initial data $(\tilde{u}_0, \tilde{u}_1)$ is actually the same 
as $(u(t_0), \partial u(t_0))$, we have an estimate for $u$
\begin{equation}
 \int \int_{\Omega_R} \frac{|u|^{p+2}}{|x|} dx dt \lesssim Q(100R)^{16(1-s)}.
 \label{estimateincone}
\end{equation}
Here $\Omega_R$ is the cone 
\[
 \Omega_R = \{(x,t): |x| + t - t_0 < R, t \in (t_0,t_0+ R)\},
\]
because in this cone $u$ and $\tilde{u}$ are the same function by the finite speed of propagation. \\
\paragraph{Lower Bound for the Integral} In the local compactness part, we have a lower bound
for the integral of $|u|^{p+2}/|x|$ in a box. Next we will show the cone $\Omega_R$ contains a lot of boxes, so we have 
a lower bound for the integral in the cone. Recall our estimate for the boxes by (\ref{boxestimate}).
\begin{equation}
 |x| + t - t_0 \leq (\frac{R_0 C_1}{d} + C_1 + 1) (t_i -t_0) < C_2 (t_i - t_0). 
\end{equation}
for all $(x,t)$ in the $i$'s box $B_i$ with $i \geq 1$. Here $C_2$ is defined to be a constant 
a little greater than $(R_0 C_1)/d + C_1 +1$. It depends only on $u$.\\
Now let us show that all boxes $B_0, B_1, B_2, \cdots, B_{m(R/C_2)}$ are completely contained in the cone 
$\Omega_R$ for each sufficiently large $R$. The cone contains the first box $B_0$ as long as $R$ is sufficiently large.
For other $1 \leq i \leq m(R/C_2)$, any point $(x,t)$ in the box $B_i$ satisfies
\[
 |x| + t - t_0 < C_2 (t_i - t_0) \leq C_2 (t_{m(R/C_2)} -t _0) \leq C_2 (R/C_2) = R.
\]
This implies this point is in the cone $\Omega_R$, by definition. Thus the cone contains those boxes we mentioned above. This gives us
\[
 \int \int_{\Omega_R} \frac{|u|^{p+2}}{|x|} dx dt \geq \sum_{i=0}^{m(R/C_2)} \int \int_{B_i} \frac{|u|^{p+2}}{|x|} dx dt. 
\]
By the local compactness result (\ref{lowbound1}) we have 
\[
\int \int_{B_i} \frac{|u|^{p+2}}{|x|} dx dt \geq \lambda(t_i)^{2 -2s} \eta_0.
\]
Thus 
\[
 \int \int_{\Omega_R} \frac{|u|^{p+2}}{|x|} dx dt \geq \sum_{i=0}^{m(R/C_2)} \lambda(t_i)^{2 -2s} \eta_0. 
\]
Now let us assume $s > 15/16$, then $2 - 2s < 1/ 8$. This means 
\[
 \lambda(t_i)^{2-2s} \geq \lambda(t_i)^{1/8},
\]
because $\lambda(t_i) \leq 1$. This implies 
\[
 \int \int_{\Omega_R} \frac{|u|^{p+2}}{|x|} dx dt \geq \sum_{i=0}^{m(R/C_2)} \lambda(t_i)^{1/8} \eta_0 = \eta_0 Q(R/C_2). 
\] 
\paragraph{Fast Growth of $Q(T)$} 
Collecting both the lower and upper bounds of the integral we have 
\[
 \eta_0 Q(R/C_2) \leq \int \int_{\Omega_R} \frac{|u|^{p+2}}{|x|} dx dt \lesssim Q(100R)^{16(1-s)}.
\]
for each sufficiently large $R > R (u, t_0)$. In other words, there exists a constant $C_u$ depending only on $u$, such that 
\[
 Q(R/C_2) \leq C_u Q(100R)^{16(1-s)}. 
\]
for large $R$. 
Let $C_3$ = $100 C_2 > 1$, we have for large $R$
\[
 Q(R) \leq C_u Q(C_3 R)^{16(1-s)}.
\]
Because $16(1-s) < 1$ when $s > 15/16$, we can choose $\kappa > 1$ such that 
\[
 16(1-s) < 1/\kappa < 1. 
\]
Using the assumption that $Q(C_3 R) \rightarrow \infty$ as $R \rightarrow \infty$, for large $R$ we have 
\[
 Q(R) \leq Q(C_3 R)^{1/\kappa}.
\]
Thus
\[
 Q(C_3 R) \geq Q(R)^\kappa. 
\]
So we have 
\[
 Q({C_3}^n R) \geq Q(R)^{\kappa^n}.
\]
Fix $R = R_2$ large so that $Q(R_2) > 1$. Then we have 
\begin{equation}
 Q({C_3}^n R_2) \geq Q(R_2)^{\kappa^n}. \label{re4}
\end{equation}
This shows that the $Q(T)$ grows very fast. This is a contradiction with the following estimate.
\paragraph{$Q(T)$ grows at a speed no faster than a linear function } By the fact
$\lambda(t) \leq 1$,
\begin{eqnarray*}
 Q(T) &=& \lambda(t_0)^{1/8} + \lambda(t_1)^{1/8} + \cdots + \lambda(t_{m(T)})^{1/8}\\
&\leq& \lambda(t_0)^{-1} + \lambda(t_1)^{-1} + \cdots + \lambda(t_{m(T)-1})^{-1} + 1.\\
& = & \frac{t_1 - t_0}{d} + \frac{t_2 - t_1}{d} + \cdots + \frac{t_{m(T)}- t_{m(T)-1}}{d} + 1\\
& = & \frac{t_{m(T)} - t_0}{d} + 1.\\
&\leq& \frac{T}{d} + 1.
\end{eqnarray*}
\paragraph{The End of the Solution} Using our linear estimate of $Q(T)$ on the left hand of 
(\ref{re4}), we have 
\[
 1 + \frac{R_2}{d} {C_3}^n \geq Q(R_2)^{\kappa^n}.
\]
for each positive integer $n$. But this is impossible for a sufficiently large $n$. This gives us 
a contradiction.

\end{document}